\documentclass[12pt]{amsart}

\usepackage[T1]{fontenc}
\usepackage[cp1250]{inputenc}
\usepackage[margin=1in]{geometry}
\usepackage{color}
\usepackage[centertags]{amsmath}
\usepackage{amsfonts}
\usepackage{amsthm}
\usepackage{newlfont}
\usepackage{amscd}
\usepackage{verbatim}
\usepackage{amsgen}
\usepackage{amssymb}
\usepackage{stmaryrd}
\usepackage{amssymb,amsmath}
\usepackage{amscd}
\usepackage{enumerate}
\usepackage{color}
\usepackage{longtable}
\usepackage{multicol}
\usepackage{xypic}
\usepackage{graphicx}
\usepackage{colordvi}
\usepackage{graphicx}
\usepackage{verbatim}
\usepackage{cancel}
\usepackage{caption}
\usepackage{epsf}
\usepackage[all,  colour, line, cmtip]{xy}
\usepackage{textcomp}

\newlength{\defbaselineskip} 
\theoremstyle{plain}
\newtheorem{thm}{Theorem}[section]
\newtheorem{cor}[thm]{Corollary}
\newtheorem{con}[thm]{Conjecture}
\newtheorem{df}[thm]{Definition}
\newtheorem{lema}[thm]{Lemma}
\newtheorem{obs}[thm]{Proposition}
\newtheorem{exm}[thm]{Example}
\newtheorem{question}[thm]{Question}

\newtheorem{fact}[thm]{Fact}

\newtheorem{rem}[thm]{Remark}
\newtheorem{pr}{Algorithm}
\theoremstyle{definition} 
\theoremstyle{definition} \newtheorem{ex}{Example}[section] %

 \numberwithin{equation}{section}

\def\z{\mathbb{Z}}
\def\n{\mathbb{N}}
\def\Z{\mathbb{Z}}
\def\N{\mathbb{N}}

\def\fa{\begin{fact}}
\def\kfa{\end{fact}}

\def\ob{\begin{obs}}
\def\kob{\end{obs}}
\def\dow{\begin{proof}}
\def\kdow{\end{proof}}

\def\tw{\begin{thm}}
\def\ktw{\end{thm}}
\def\hip{\begin{con}}
\def\khip{\end{con}}
\def\lem{\begin{lema}}
\def\klem{\end{lema}}
\def\ex{\begin{exm}}
\def\prog{\begin{pr}}
\def\kprog{\end{pr}}
\def\wn{\begin{cor}}
\def\kwn{\end{cor}}
\def\uwa{\begin{rem}}
\def\kuwa{\end{rem}}
\def\kex{\end{exm}}
\def\dfi{\begin{df}}
\def\kdfi{\end{df}}
\def\Gg{\mathcal G}
\def\Gf{\mathcal G}
\definecolor{zielony}{rgb}{0.5, 0.9, 0.1}
\definecolor{czerwony}{rgb}{0.9, 0.2, 0.1}
\definecolor{niebieski}{rgb}{0.3, 0.1, 0.9}

\begin{document}
\title{Finite phylogenetic complexity of $\Z_p$ and invariants for $\z_3$}
\author{Mateusz Micha\l ek}
\thanks{Supported by a grant Iuventus Plus of the Polish Ministry of Science, project 0301/IP3/2015/73.}
\begin{abstract}
We study phylogenetic complexity of finite abelian groups - an invariant introduced by Sturmfels and Sullivant \cite{SS}. The invariant is hard to compute - so far it was only known for $\z_2$, in which case it equals $2$ \cite{SS, Sonja}. We prove that phylogenetic complexity of any group $\z_p$, where $p$ is prime, is finite. We also show, as conjectured by Sturmfels and Sullivant, that the phylogenetic complexity of $\z_3$ equals $3$.
\end{abstract}
\maketitle
\section{Introduction}
The motivation for our work comes from phylogenetics - a science that aims at reconstructing the history of evolution. We will not present here all the concepts from phylogenetics as they are not needed for the statement and the solution of the problem that we study. Let us just say that to any tree $T$ and a finite abelian group $G$, by considering a Markov process on a tree, one associates a projective toric variety $X(T,G)$. The explicit description of the variety and the associated polytope is given in Definition \ref{def:poly}. We refer the interested reader to \cite{pachsturm, 4aut, SS, jaPhD, jajalg}, where the relations to phylogenetics and applications are explained in detail. The equations defining $X(T,G)$ are called phylogenetic invariants. In all the cases that we study, determining phylogenetic invariants for any tree $T$ was reduced to so-called star or claw trees using toric fiber product \cite[Theorem 26]{SS}, \cite[Corollary 2.11]{Sethtfp}. These trees, denoted by $K_{1,n}$ have one inner vertex and $n$ leaves. Let us cite Draisma and Kuttler \cite{DK}:

"We have now reduced the ideals of our equivariant models to those
for stars, and argued their relevance for statistical applications. The
main missing ingredients for successful applications are equations for
star models. These are very hard to come by (...)".

In our previous work with Maria Donten-Bury \cite{DBM} we have shown how to obtain phylogenetic invariants of bounded degree. However, it is highly nontrivial to obtain such a bound. To study these bounds Sturmfels and Sullivant defined two functions.
\dfi[$\psi(n,G), \psi (G)$]
Let $\psi(n,G)$ be the degree in which the (saturated) ideal defining $X(K_{1,n},G)$ is generated. Let $\psi(G)$, called the phylogenetic complexity of $G$, be the supremum of $\psi(n,G)$ over $n\in\n$.
\kdfi
As observed by Sturmfels and Sullivant \cite{SS}: "The phylogenetic complexity $\psi(G)$ is an intrinsic invariant of the group $G$. (...) It would be interesting to study the group-theoretic meaning of this invariant."
However, these invariants are very hard to compute. So far we only know $\phi(\z_2)=2$ \cite{SS, Sonja}. Based on numerical computations Sturmfels and Sullivant proposed the following conjecture.
\hip\cite[Conjecture 29]{SS}\label{conSS}
For any finite abelian group $G$ we have $\psi(G) \leq |G|$.
\khip
However, for $G\neq \z_2$ we do not know if $\psi(G)$ is finite. Our first main theorem is as follows.
\tw\label{thm:zp}
For any prime number $p$
the phylogenetic complexity of $\z_p$ is finite.
\ktw
Depending how general the model is there are other qualitative results on the degree of phylogenetic invariants. For very general, so-called equivariant models, the fact that on set-theoretic level there exists a bound was proved in \cite{DK, DK2, DE}. For the class of $G$-models that includes all the models introduced in this article, on the level of projective schemes the bounds were obtained in \cite{JaJCTA}. Finally, for group-based models, but \emph{only on Zariski open set,} the bound of the degrees by $|G|$ was proved in \cite{CFSM}. Our second main theorem is as follows.
\tw\label{thm:z3}
The phylogenetic complexity of the group $\z_3$ equals $3.$
\ktw
This allows to find all phylogenetic invariants for any tree for the group $\z_3$. As far as we know, this is the only model, different from the Jukes-Cantor model, where the complete list of phylogenetic invariants for any tree is obtained. For real data applications of phylogenetic invariants we refer for example to \cite{RusinkoHipp}. We would also like to mention that a related result was recently obtained by Donten-Bury in \cite{Marysianowyz3} on scheme-theoretic level.

The techniques that we use rely entirely on algebraic combinatorics. We present the above described problems in the combinatorial terms in Section \ref{Definitions}. In different words, we study algebraic properties of a family of integral polytopes.

Although the original construction of varieties $X(T,G)$ was inspired by phylogenetics, recently they appeared in other sciences \cite{Man, Man2, Man3, SX}. We would like also to mention that the varieties $X(T,G)$ share many other very interesting algebraic and combinatorial properties related to their Hilbert polynomial, normality and deformations \cite{BBKM, BW, kaieK, Rusinko}.

The problems of the degrees in which toric ideals are generated appear in many different contexts \cite{Bruns}. Let us summarize the results and conjectures about group-based models in the following table. 
\newline
\begin{tabular}{|p{3,1 cm}|p{1,8 cm}|p{2,3 cm}|p{2,3 cm}|p{2,25 cm}|p{2,8 cm}|}
\hline
&\multicolumn{5}{ |c| }{Group-based Models} \\
\hline
polynomials defining:& $\z_2$ & $\z_3$ & $\z_2\times\z_2$ & $\z_p$ & $G$ \\
\hline
Gr\"obner basis & degree $2$ by \cite{Sonja} &&&& Question \ref{Question1}\\
\hline
generators of the ideal & degree $2$ by \cite{SS} & degree $3$ by Theorem \ref{thm:z3} & Conjecture \cite[Conjecture 30]{SS} & finite by Theorem \ref{thm:zp} & Conjecture \ref{conSS} \cite[Conjecture 29]{SS} \\
\hline
the projective scheme & & degree $3$ \cite{Marysianowyz3} & degree $4$ \cite{JaJCTA} & & finite by \cite{JaJCTA}\\
\hline
set-theoretically&&&&&finite by \cite{DE}
\\
\hline
on a Zariski open subset & && degree $4$ \cite{JaAdvGeom} && degree $\leq |G|$ \cite{CFSM}\\
\hline
\end{tabular}

As one can see the higher the row, the finer algebraic properties are required. On the other hand columns to the right provide bigger and more general groups. This provides our table with "diagonal" structure with theorems mostly on and below the diagonal and conjectures above it. Taking this into account the following question, for now out of reach, is the most difficult.
\begin{question}\label{Question1}
What are the bounds on the degree of Gr\"obner basis for group-based models?
\end{question}

\section*{Acknowledgements}
The project was started when the author was visiting Freie Universit\"at in Berlin thanks to Mobilno\'s\'c Plus program of the Polish Ministry of Science. The author is also a member of AGATES group and PRIME DAAD fellow. We would like to thank two reviewers for suggesting many improvements.

\section{Definitions}\label{Definitions}
Throughout the article $G$ will be a finite abelian group. Let $[n]=\{1,\dots,n\}$.
\dfi[Flow]
Fix $n\in \n$. A group-based flow $f$ (on $n$) or simply a flow is a function $f:[n]\rightarrow G$ represented by an $n$-tuple of group elements $(f(1)=g_1,\dots,f(n)=g_n)$ and satisfying $\sum_{i=1}^n f(i)=0\in G$. With a coordinate-wise action flows form a group of flows $\Gf$ isomorphic to $G^{n-1}$.

We will say that an element $g\in G$ belongs to a flow $f$ if it belongs to the image of $f$. We let $\underline 0=(0,\dots,0)$ be the neutral element of $\Gf$.

\kdfi
The object of our study is a family of integral polytopes indexed by an integer $n\in \n$ and a finite abelian group $G$. These polytopes are combinatorial objects representing the group of flows $\Gf$. They are subpolytopes of a unit cube, hence all their integral points are vertices. The vertices are in bijection with elements of $\Gf$.
\dfi[Polytope $P_{n,G}$]\label{def:poly}
Consider the lattice $M\cong \z^{|G|}$ with a basis corresponding to elements of $G$. Consider $M^n$ with the basis $e_{(i,g)}$ indexed by pairs $(i,g)\in [n]\times G$. We define an injective map of sets:
$$\Gf\rightarrow M^n,$$
by $f\rightarrow \sum_{i=1}^n e_{(i,f(i))}$. The image of this map defines the vertices of the polytope $P_{n,G}$.
\kdfi
\ex\label{ex:z2}
For $n=3$ and $G=\z_2$ we have four flows:
$$(0,0,0),(0,1,1),(1,0,1),(1,1,0)\in \Z_2\times \Z_2\times \Z_2.$$

Hence, the polytope $P_{3,\z_2}$ has the following four vertices corresponding to the flows above:
$$(1,0,1,0,1,0),(1,0,0,1,0,1),(0,1,1,0,0,1),(0,1,0,1,1,0)\in\Z^2\times \Z^2\times \Z^2,$$
where $(1,0)\in \Z^2$ corresponds to $0\in \Z_2$ and $(0,1)\in\Z^2$ corresponds to $1\in \Z_2$.
\kex
The operation of addition is \emph{different} for the two representations. Under the first representation flows form the finite group $\Gf$ -- in particular, every element has finite order (equal to two in Example \ref{ex:z2}). Under the second representation the addition is simply addition of integer vectors induced from $\z^{n|G|}$. In particular, every element has infinite order (and in Example \ref{ex:z2} the four vertices are independent).

The polytope $P_{n,G}$ does not have to be normal (however for many groups it is conjectured to be). The associated toric variety in the sense of \cite{Stks} is isomorphic to $X(K_{1,n},G)$ and is the main object of our study. However, as we will see, the language of flows, due to the group structure, is easier and will be used throughout the article instead of the language of polytopes and integral points. We refer the reader to \cite{Fult, CLS} for background on toric varieties.

Phylogenetic invariants, that is equations defining $X(K_{1,n},G)$, correspond to integral relations among vertices of $P_{n,G}$. We use the following combinatorial restatement. We say that two multisets of flows (on $n$) $M_1$ and $M_2$ are \emph{compatible} if for any $i\in [n]$ we have $\cup_{f\in M_1} \{f(i)\}=\cup_{g\in M_2} \{g(i)\}$ as multisets of elements of $G$. This is equivalent to the fact that corresponding vertices sum up to the same lattice element. Of course to two compatible multisets we can add a flow, obtaining two bigger compatible multisets. The degree of a binomial corresponding to compatible multisets $M_1$, $M_2$ equals the cardinality of any of the multisets (note that both have to be of the same cardinality). By \emph{exchanging} a multiset of group based flows we will always mean exchanging it with a compatible multiset.

\ex\label{ex:comp}
For $G=\z_2$ and $n=6$ using the representation of flows as tuples of group elements we have e.g.~the following two compatible multisets:
$$M_1=((1,1,1,1,1,1),(0,0,0,0,0,0), (1,1,1,1,0,0))$$
and
$$M_2=((0,1,0,1,0,0),(1,1,1,0,1,0),(1,0,1,1,0,1)). $$
We could represent them as tables:
$$\xymatrix{T_1&=&\txt{{(1,1,1,1,1,1)}\\{(0,0,0,0,0,0)}\\{(1,1,1,1,0,0)}}}$$
and
$$\xymatrix{T_2&=&\txt{{(0,1,0,1,0,0)}\\{(1,1,1,0,1,0)}\\{(1,0,1,1,0,1)}}}.$$
Note that \emph{any} two flows from $M_1$ are \emph{not compatible} with any two flows from $M_2$. However, $((1,1,1,1,1,1),(0,0,0,0,0,0))$ is compatible with $((0,1,0,1,0,0),(1,0,1,0,1,1))$, hence we may \emph{exchange} them obtaining:
$$\tilde M_1=((0,1,0,1,0,0),(1,0,1,0,1,1),(1,1,1,1,0,0)).$$
Now $M_1$ and $\tilde M_1$ are \emph{compatible} and the last two flows in $\tilde M_1$ are compatible with the last two flows of $M_2$. Hence, we have a sequence of compatible multisets $M_1\sim \tilde M_1\sim M_2$.

We started from a degree three binomial and generated it using degree two binomials. 
\kex

The ideal of the variety $X(K_{1,n},G)$ is generated in degree $d$ if and only if the compatiblilty relation on multisets equals the transitive closure of the restriction of compatibility relation to multisets of cardinality at most $d$ and the operation of adding a flow. More explicitly, if and only if we are able to pass from any multiset to any compatible multiset in a series of steps, each time exchanging a submultiset (of one multiset) with at most $d$ flows by a compatible multiset of flows. By exchanging two flows $f,g$ (in one multiset) on a set of indices $I$, we will mean replacing them with two flows $f',g'$ such that $f'(i)=f(i)$ and $g'(i)=g(i)$ for $i\not\in I$, $f'(i)=g(i)$ and $g'(i)=f(i)$ for $i\in I$. Notice that this is only possible if $\sum_{i\in I} f(i)=\sum_{i\in I} g(i)$.
\ex
Using the notation as in Example \ref{ex:comp} while passing from $M_1$ to $\tilde M_1$ we exchanged $f=(1,1,1,1,1,1)$ and $g=(0,0,0,0,0,0)$ on indices $I=\{1,3,5,6\}$ obtaining $f'=(0,1,0,1,0,0)$ and $g'=(1,0,1,0,1,1)$.

\kex
\section{Bounded phylogenetic complexity for $\z_p$}
We hope that the arguments of this section will be generalized to arbitrary finite abelian groups in future work.
The whole section is devoted to the proof of Theorem \ref{thm:zp}. 
We will prove an equivalent statement: for $n$ large enough, the ideal corresponding to the claw tree with $n+1$ leaves is generated in the same degree as for $n$ leaves, that is $\psi(n+1,\z_p)=\psi(n,\z_p)$.

Consider two compatible multisets $M_1$ and $M_2$ of flows on $n+1$. The proof that one can pass from $M_1$ to $M_2$ exchanging at most $\psi(n,\z_p)$ flows at each step is inductive on the cardinality of $M_1$ (that is the same as the cardinality of $M_2$). The case when $M_1$ is of cardinality one (or at most $\psi(n,\z_p)$) is trivial.

Otherwise, choose $f_1\in M_1$ and $g_1\in M_2$. Suppose that $f_1$ and $g_1$ agree on $k$ indices. In our proof we inductively increase $k$. If $k=n+1$ we can conclude (as $f_1=g_1$) by reducing to the case of smaller cardinality of the multisets.

We distinguish two basic cases. Briefly, in the first one the flows $f_1$ and $g_1$ differ a lot, i.e.~$k$ is small. Here we can easily conclude using Lemma \ref{lem:zamiana} that will also be useful in the following Section \ref{sec:z3}. In the second case we first consider the situation in which all the flows in $M_1$ and $M_2$ are very much alike. This is similar e.g.~to the approach presented in \cite{JaMichalMatroid}. The application of the results of \cite{JaJCTA} allows us to finish the proof.

\lem\label{lem:zamiana}
Fix $G=\z_p$ for $p$ a prime number. Suppose there exist two flows $f,g$ and a set of indices $I$ of cardinality $p-1$ such that for each $i\in I$ we have $f(i)\neq g(i)$. Then for any subset of indices $I'$ disjoint from $I$ there exists a subset $I''\subseteq I$ such that we can exchange $f$ and $g$ on $I'\cup I''$.
\klem
\dow
One possible proof, that we leave as an exercise for the reader, is by induction. We  propose a different one based on Combinatorial Nullstellensatz \cite{Alon, JaNull, Michal}.


Taking $-h:=\sum_{i\in I'} f(i)-g(i)$, it is enough to show that 
there exists a subset $I''\subset I$ such that $\sum_{i\in I''}f(i)-g(i)=h$. 
Let us consider the polynomial:
$$\prod_{h'\neq h}((\sum_{i\in I} (f(i)-g(i))x_i)-h')\in \z_p[x_i]$$
with nonzero coefficient of $\prod_{i\in I}x_i$. By Alon's Nullstellensatz there exists a point $P\in\{0,1\}^{|I|}$ on which the polynomial does not vanish. The coordinates of $P$ which are nonzero identify a subset $I''\subset I$ with the desired property.
\kdow


By the action of $\Gg$ (i.e.~adding a fixed flow to all other flows) we may assume $g_1=\underline 0$. Let $I:=\{i:f_1(i)\neq 0\}$. For each $i\in I$ there must exist a flow $f\in M_1$ such that $f(i)=0$.

 We first conclude in the easy case:

Case 1) $|I|>3p^2$.

As all nonzero elements of $G$ are indistinguishable, we may assume that there are $3p$ indices $i$, such that $f_1(i)=1$. Let $I'$ be the set of such indices. For any $j\in I'$ we may assume that any flow $f\in M_1$ such that $f(j)=0$ has to assign one to all but at most $p-1$ indices from $I'$, i.e.~it has to agree with $f_1$ on $I'$ on all apart from at most $p-1$ indices. Indeed, otherwise by Lemma \ref{lem:zamiana}, we can exchange $f$ and $f_1$ on a subset of $I'$ that contains $j$ and increase $k$. Also other flows $f'\in M_1$ must assign $1$ to at least $3p-2p=p$ indices of $I'$ as otherwise we would be able to exchange them first with an $f$ such that $0\in f(I')$ and then conclude as above. As $M_1$ and $M_2$ are compatible, there must exist a flow $g\in M_2$ that assigns one to at least $p$ elements in $I'$. We can make an exchange between $g$ and $g_1$ increasing $k$.

 Hence, from now on we assume:
 
Case 2) $|I|\leq 3p^2$.

Our aim is to reduce to the situation, in which there exists one flow $f_0\in M_1$ such that $f_0(i)=0$ for all $i\in I$ (or similarly a flow $g_0\in M_2$ such that $g_0(i)=f_1(i)$ for all $i\in I$). This will finish the proof, as then we can exchange $f_1$ and $f_0$ on $I$ (recall that $f_1$ is a flow and $f_1(s)=0$ for $s\not\in I$). Hence, we fix a flow $f_0$.
Our strategy in this case is to apply results of \cite{JaJCTA} to prove the main Lemma \ref{lem:2flows}. Once it is proved, one can easily increase $|\{i\in I: f_0(i)=0\}|$.
It is often useful to encode flows as colorings.

\dfi[Coloring]
A coloring of length $n$ is a function $f:[n]\rightarrow[g]\cup \{0\}$. The number $g$ is called the number of colors. The support of the coloring $f$ is defined as $\{k\in[n]: f(k)\neq 0\}$.
\kdfi
\dfi[Transformation]
Consider two colorings $f_1,f_2:[n]\rightarrow[g]\cup\{0\}$. Suppose that there exist two numbers $0\leq k_1,k_2\leq n$ such that $k_j$ is \emph{not} in the support of $f_j$. Moreover suppose $f_1(k_2)=f_2(k_1)$. Define $f_j'(x)=f_j(x)$ for $x\neq k_1,k_2$. Moreover, $f_1'(k_j):=f_2(k_j)$ and $f_2'(k_j):=f_1(k_j)$. We call $f_1',f_2'$ a \emph{transformation} of $f_1$ and $f_2$.
\kdfi
Transformation of colorings corresponds to exchanging the fixed color in two colorings with the $0$ color. A multiset of colorings can be transformed into another by choosing two colorings and transforming them. We generate an equivalence relation on multisets of colorings by transformations. Abusing the notation the relation is also called transformation. We note that transformations of colorings give rise to compatible exchanges of multisets of flows (generated by quadrics).
\lem\label{combinatorial}\cite[Lemma 6.5]{JaJCTA}\label{lem:JCTA}
Let us fix three natural numbers: $g$ (number of colors), $s$ (bound on the support) and $a\geq 2$. Fix $\epsilon >0$. There exists $N\in \n$, such that for all $n\geq N$ any collection of colorings $f_1,\dots,f_m:[n]\rightarrow [g]\cup\{0\}$ with support of cardinality at most $s$ can be transformed into a collection $f_1',\dots,f_m'$ with the following property:

 there exist $\lfloor(1-\epsilon)\frac{n}{a}\rfloor$ numbers $x<n$ divisible by $a$, such that for any $f_j'$ and any $x$ at most one of the numbers $x,x+1,\dots,x+a-1$ is in the support of $f_j'$.
\klem

\lem\label{lem:2flows}
We may assume that there exist two flows $f_3,f_4\in M_1$ different from $f_1$, such that $|j\not\in I:f_3(j)\neq f_4(j)|>2p$.
\klem
\dow
Let us assume the contrary. By exchanging $M_1$ and $M_2$ we may assume that after restricting all flows to the complement of $I$:
\begin{enumerate}
\item any two flows $f_i,f_j$ different from $f_1$ differ on at most $2p$ indices,
\item any two flows $g_i,g_j$ different from $g_1$ differ on at most $2p$ indices.
\end{enumerate}
We will show that in this case we can transform $M_1$ to $M_2$ exchanging at most $\phi(n,G)$ flows at a time.
We proceed in the following steps:
\begin{enumerate}[1)]
\item Show that there exists a function $q:[n+1]\setminus I\rightarrow G$ such that each flow $f\in M_1\setminus f_1$ \emph{and} each flow $g\in M_2\setminus g_1$ differ from $q$ (on the complement of $I$) on at most $4p$ indices.
\item Apply Lemma \ref{lem:JCTA}, to obtain two indices $i_0,i_1$, such that each flow $f\in M_1\setminus f_1$ \emph{and} each flow $g\in M_2\setminus g_1$ differ from $q$ on $i_0,i_1$ on at most $1$ index.
\item Reduce to the case of flows on $n$.
\item Apply induction by lifting the exchanges among flows on $n$ and generate the relation among flows on $n+1$.
\end{enumerate}

{\bf Step 1}:

Let $q$ be the restriction of any flow in $M_1\setminus f_1$ to the complement of $I$.
By assumption $(i)$ any other flow in $M_1\setminus f_1$ differs from $q$ on at most $2p<4p$ indices. We want to show that also the flows in $M_2\setminus g_1$ cannot differ from $q$ on more than $4p$ indices.

We act (coordinate-wise) on $M_1\setminus f_1$ with the inverse of the flow that defined $q$ and restrict to the complement of $I$. Identifying $0\in G$ with $0\in\N$ and other group elements with consecutive natural numbers
we obtain a multiset of colorings with the support bounded by $2p$. It remains to show that after this action also the elements of $M_2\setminus g_1$ must have a bounded support. Indeed, on average every element of $M_2\setminus g_1$ must have at most $2p$ nonzero elements (notice that after restricting to the complement of $I$ the multiset $M_1\setminus f_1$ is compatible with $M_1\setminus g_1$). However, as the flows differ on at most $2p$ indices all flows must have at most $4p$ nonzero elements.

{\bf Step 2}:

Taking $s=4p$, $a=2$, $\epsilon<1/2$ and applying Lemma \ref{lem:JCTA} to flows/colorings from $M_1\setminus f_1$ and $M_2\setminus g_1$ restricted to the complement of $I$ we may assume that there exist:
\begin{enumerate}
\item indices $i_0,i_1\not\in I$,
\item two group elements $g:=q(i_0),h:=q(i_1)$,
\end{enumerate}
 such that
 \begin{enumerate}
\item any flow $f_j$ for $j\neq 1$ satisfies $f_j(i_0)=g$ or $f_j(i_1)=h$,
\item any flow $g_j$ for $j\neq 1$ satisfies $g_j(i_0)=g$ or $g_j(i_1)=h$.
\end{enumerate}

{\bf Step 3}:

By permuting columns, for simplicity of notation assume $i_1=i_0+1$.
We may replace the flows $f_j,g_j$ by flows $f_j',g_j'$ on $n$ where:
$$f'_j(k)=
\begin{cases}
f_j(k)\quad k<i_0\\
f_j(k+1)\quad k>i_0
\end{cases},
f'_j(i_0)=f_j(i_0)+f_j(i_1)$$$$
g'_j(k)=
\begin{cases}
g_j(k)\quad k<i_0\\
g_j(k+1)\quad k>i_0
\end{cases},
g'_j(i_0)=g_j(i_0)+g_j(i_1).
$$
The obtained multisets remain \emph{compatible}. Indeed, restricting the flows to indices $i_0,i_1$  we obtain exactly the same pairs of group elements for $M_1$ and $M_2$ (these pairs are either $(g,h)$ or $(g,h')$ for $h'\neq h$ or $(g',h)$ for $g'\neq g$).

{\bf Step 4}:

By induction, for flows on $n$ we only need to exchange at most $\phi(n,G)$ flows at each step. Notice, that each exchange lifts to an exchange among $M_1$ and $M_2$. Indeed, the exchanges on the summed entries lift to exchanges among pairs of indices $i_0, i_1$. Notice that this is \emph{not} enough to conclude, as at the end we do not obtain the same multisets - just the multisets that after summing up entries $i_0$ and $i_1$ are the same. However, the entries on $i_0$ and $i_1$ may be adjusted using simple quadratic moves (i.e.~exchanges among two flows), which finishes the proof of the lemma.
\kdow
We may now increase the number of indices $i\in I$ such that $f_0(i)=0$. Choose $i_0\in I$ such that $f_0(i_0)\neq 0$ and $f'\in M_1$ such that $f'(i_0)=0$. Exchanging $f_0$ either with $f_3$ or $f_4$ from the lemma we may assume that $f_0$ and $f'$ differ on at least $p-1$ indices not from $I$. By Lemma \ref{lem:zamiana} we can exchange $f_0$ and $f'$ on index $i$ and some indices not in $I$. This finishes the proof of Theorem \ref{thm:zp}.
\section{The phylogenetic complexity of $\z_3$}\label{sec:z3}
The whole section is devoted to the (very technical) proof of Theorem \ref{thm:z3}.
The main methods are as follows. For small $n$ we use direct computational results using \cite{Normaliz, 4ti2, DBM}. For large $n$ the basic distinction into the case where two flows are very different or very much alike remains valid. However, in both cases (especially when two flows are much alike) we failed to find short, easy proofs or conclude simply by deeper and deeper case study. Our arguments rely on a mixture of counting arguments and combinatorial tricks. In particular, we believe that it is impossible to follow this proof without a pen and a sheet of paper. The author would much appreciate an approach that would be significantly simpler (as this would rise hope for attacking larger groups, e.g.~\cite[Conjecture 30]{SS}).

Our proof is inductive on $n$, the length of the flows. For fixed $n$ we show inductively that any binomial of degree $d>3$ is generated by cubics. For $n<7$ it was known before that the ideal of $\z_3$ is generated by cubics \cite{DBM}. Thus we assume $n\geq 7$.
Let us fix two compatible multisets of group-based flows $M_1=\{f_1,\dots,f_d\}$ and $M_2=\{g_1,\dots,g_d\}$. To simplify the notation, by the action of $\Gf$ we may assume that $g_1=\underline{0}$ is the trivial flow. Suppose that $f_1$ and $g_1$ agree on $k-1$ indices. Our aim is to increase $k$, until $k-1=n$. Of course if $k-1=n$ then $g_1=f_1$ and we may conclude by the induction on $d$. By the action of the symmetric group $S_n$ we may assume that $f_1(i)=0$ if and only if $1\leq i\leq k-1$. As the proof is quite complicated we decided to include diagrams that describes most important cases. While reading the proof we encourage the reader to follow at which node we are. The proofs are "depth-first, left-first".

\subsection{Suppose $k\neq n-1$.}\label{sec:3}

There are two possibilities we consider. Either $f_1$ attains the same element for at least three indices greater or equal to $k$ or $k=n-3$ and we can assume $f_1(k)=f_1(k+1)=1$, $f_1(k+2)=f_1(k+3)=2$.
\subsubsection{Suppose (by the action of the nontrivial automorphism of $\z_3$) we have:\\ $f_1(k)=f_1(k+1)=f_1(k+2)=1$.}\label{subsec:1}
The diagram for this subsection is presented below. In each vertex we presented main case assumptions and lemmas with very short comments that the reader may unravel while following the proof. 
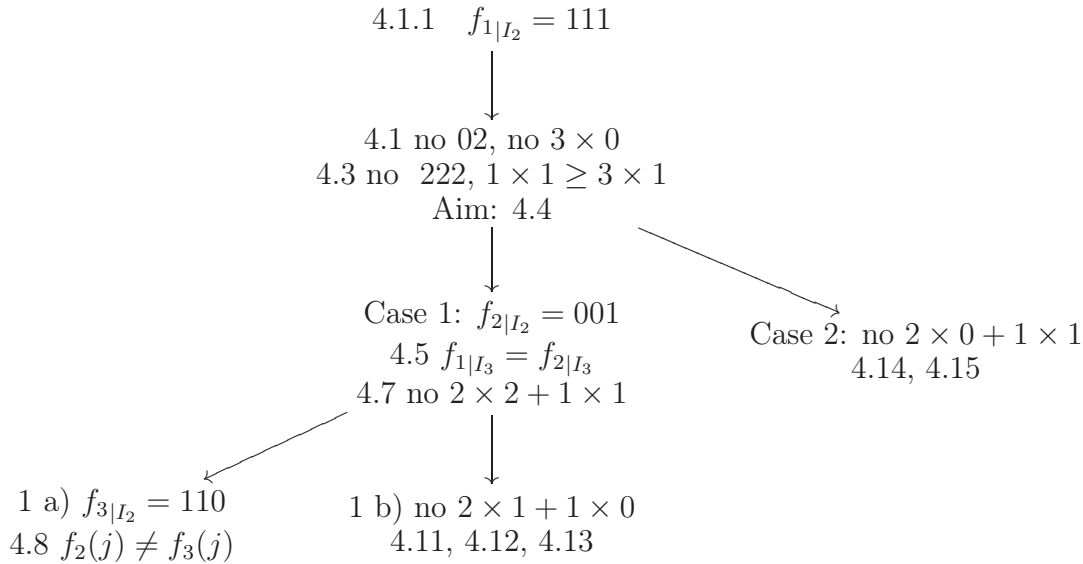
\[
\xymatrix{
&\text{\ref{subsec:1}} \quad {f_1}_{|I_2}=111\ar[d]&&\\
&{\txt{{\text{\ref{lem:no02} no }02,\text{ no }$3\times 0$}\\ {\text{\ref{lem:more1} no } 222, $1\times 1\geq3\times 1$}\\{\text{Aim: }\ref{ob:gl}} }}\ar[d]\ar[rd]&&\\
&\txt{{Case 1: ${f_2}_{|I_2}=001$}\\{\ref{lem:f_2zgf_1} ${f_1}_{|I_3}={f_2}_{|I_3}$}\\{\ref{lem:no221} no $2\times 2+1\times 1$}}\ar[ld]\ar[d]&\txt{{Case 2: no $2\times 0+1\times 1$}\\{\ref{lem:all}, \ref{lem:two}}}&\\
\txt{{1  a) ${f_3}_{|I_2}=110$}\\{\ref{lem:f2vsf3}  $f_2(j)\neq f_3(j)$}}&\txt{{1 b) no $2\times 1+1\times 0$}\\{\ref{lem:1agree}, \ref{lem:0agree}, \ref{lem:2agree}}}&\\
}
\]

Let us group the indices:
$$I_1:=\{1,\dots,k-1\}, I_2:=\{k,k+1,k+2\}, I_3:=\{k+3,\dots,n\}.$$
Note that $I_1$ or $I_3$ may be empty.
\lem\label{lem:no02}
If there exists a flow $f\in M_1$ that has either
\begin{enumerate}
\item  both $0$ and $2$ or
\item only $0$,
\end{enumerate}
 in the image of $I_2$, then we can increase $k$.
\klem
\dow
In the first case, we may exchange the elements on the preimages of $0$ and $2$ in $I_2$ between $f$ and $f_1$, increasing the number of neutral elements in $f_1$.
In the second case we exchange $f$ and $f_1$ on $I_2$.
\kdow
\uwa\label{rem:001l011}
As $g_1$ has only neutral elements in the image of $I_2$, for each $i\in I_2$ there must exist a flow $f\in M_1$, such that $f(i)=0$. By Lemma \ref{lem:no02} such flows $f$ can be only of two types; either they have precisely twice $1$ and once $0$ in the image of $I_2$ or twice $0$ and once $1$.
\kuwa
\lem\label{lem:more1}
We may assume that:
\begin{enumerate}
\item there is no flow among $M_1$ that associates three times $2$ on $I_2$,
\item there exist at least as many flows in $M_1$ that associate to $I_2$ exactly one $1$ as those that associate three times $1$.
\end{enumerate}
\klem
\dow
Suppose that $f\in M_1$ associates three times $2$ on $I_2$. There exists a flow $f'\in M_1$ that associates either twice $1$ and once $0$ or twice $0$ and once $1$ on $I_2$. We can make an exchange between $f$ and $f'$ on two indices from $I_2$ on which $f'$ has value $0$ and $1$. This proves the first part of the lemma by Lemma \ref{lem:no02}.

Hence, also by Lemma \ref{lem:no02} we see that each flow in $M_1$ has at least one $1$ on $I_2$. If there are strictly more flows that associate three times $1$ than those that associate $1$ only once (with each of the other two entries equal either to $0$ or $2$), then an average of associations of $1$ on $I_2$ for the whole $M_1$ is above $2$. As $M_1$ and $M_2$ are compatible the same must be true for $M_2$. Hence, there must exist $g\in M_2$ that associates to $I_2$ three times $1$. We can make an exchange between $g$ and $g_1$ on $I_2$, increasing $k$ by three.
\kdow
Our aim is to prove the following proposition.
\ob\label{ob:gl}
We can compatibly change $M_1$, so that
there exists a function $h:I_1\cup I_3\rightarrow\z_3$ such that:
\begin{enumerate}
\item for every index $j_0\in I_3$ we have $h(j_0)= f_1(j_0)$,
\item for any $f\in M_1\setminus f_1$ there exists at most one index $j_0\in I_1\cup I_3$ such that $h(j_0)\neq f(j_0)$,
\end{enumerate}
as otherwise we are able to increase $k$.
\kob
The proof of the above crucial proposition is divided into several parts.


{\bf Case 1)} There exists a flow in $M_1$, say $f_2$, that in the image of $I_2$ has precisely two $0$ and one $1$.
By the action of $S_{I_2}\subset S_n$ we may assume $f_2(k)=f_2(k+1)=0$ and $f_2(k+2)=1$.
In this case the function $h$ in Proposition \ref{ob:gl} will be the restriction of $f_2$ to $I_1\cup I_3$.

\lem\label{lem:f_2zgf_1}
For any $i\in I_3$ we can assume $f_1(i)=f_2(i)$. For any $i\in I_1$ we have either $f_2(i)=0$ or $f_2(i)=1$. The number of $i\in I_1$ such that $f_2(i)=1$ equals $2$ modulo $3$.
\klem
\dow
Suppose there is $i\in I_3$ such that $f_1(i)\neq f_2(i)$. We may exchange $f_1$ and $f_2$ on $i$ and on one or two indices $k, k+1$, so that we get two flows. This would increase $k$.

If there is $i\in I_1$ such that $f_2(i)=2$ then we can exchange $f_1$ and $f_2$ on $i,k,k+1$. This would also increase $k$. The last statement follows from the first, because $f_2$ is a flow.
\kdow
\wn
Point 1) of Proposition \ref{ob:gl} holds.\hfill$\square$
\kwn
By Remark \ref{rem:001l011} there must exist a flow in $M_1$, say $f_3$, such that $f_3(k+2)=0$. 
\lem\label{lem:no221}
We may assume that there are no flows in $M_1$ that to $I_2$ associate twice $2$ and once $1$. In particular, by Lemma \ref{lem:no02}, all flows that associate only one $1$ on $I_2$ associate exactly two $0$.
\klem
\dow 
Suppose such a flow $f$ exists. We may assume that $f(k)=f(k+1)=2$ and $f(k+2)=1$ as otherwise we could conclude by Lemma \ref{lem:no02}, by exchanging $f$ and $f_2$ on $k$, $k+1$. As above we must have $f_3(k)=f_3(k+1)=1$. Exchanging $f_3$ and $f_2$ on $k+1$ and $k+2$ and then exchanging $f$ with $f_2$ on $k$ and $k+2$ we can conclude by Lemma \ref{lem:no02}. 
\kdow

By Remark \ref{rem:001l011} we may consider the following two cases.

Case 1 a) $f_3(k)=f_3(k+1)=1$. 
\lem\label{lem:f2vsf3}
On $I_1\cup I_3$ any two flows $f_{i_1}$ and $f_{i_2}$ such that $\sum_{i\in I_2}f_{i_1}(i)=1$ (e.g.~$f_2$) and $\sum_{i\in I_2}f_{i_2}(i)=2$ (e.g.~$f_3$) differ on $I_1\cup I_3$ on precisely one index $j$ for which $f_{i_1}(j)=f_{i_2}(j)-2$.
\klem
\dow
By making exchanges with $f_2$ and $f_3$ it is enough to prove the lemma for these two flows.
Suppose there exists an index $i\in I_1\cup I_3$ such that $f_{2}(i)=f_{3}(i)+2$. We may then exchange the flows $f_2$ and $f_3$ on the indices $i$ and $k+2$. We obtain a flow that associates $0$ to all the indices of $I_2$, which allows to increase $k$ by Lemma \ref{lem:no02}.

Suppose there is more than one index $j_1, j_2\in I_1\cup I_3$ such that $f_2(j_i)=f_3(j_i)+1$ for $i=1,2$. Then, as above, we may exchange $f_2$ and $f_3$ on $j_1,j_2$ and $k+2$. This finishes the proof - one such index must exist as $f_2$ and $f_3$ are flows.
\kdow
\dfi[index $j$]
We define the distinguished index $j$ to be the unique $j\in I_1\cup I_3$ such that $f_2(j)\neq f_3(j)$.
\kdfi

We now prove point 2) of Proposition \ref{ob:gl}. Consider any flow $f\in M_1$ different from $f_1,f_2,f_3$.

If $\sum_{i=k}^{k+2}f(i)=2$ then by Lemma \ref{lem:f2vsf3} we see that $f$ differs from $f_2$ on $I_1\cup I_3$ on precisely one index, that allows to conclude in this subcase.

If $\sum_{i=k}^{k+2}f(i)=1$ we may exchange $f$ and $f_2$ on $I_2$. By Lemma \ref{lem:f2vsf3} we see that $f$ differs from $f_3$ on $I_1\cup I_3$ on precisely one index. As $f_2$ and $f_3$ agree on $I_1\cup I_3\setminus\{j\}$, we have two possibilities. Either $f$ agrees with $f_2$ on $I_1\cup I_3$ or there exists $j'\in I_1\cup I_3\setminus\{j\}$ such that $f(j)=f_2(j)+2$ and $f(j')=f_2(j')-2$. We have to exclude the latter case. By Lemma \ref{lem:no02} we have $f(I_2)=\{1,1,2\}$ or $f(I_2)=\{0,0,1\}$. By the $S_{I_2}$ action it is enough to consider the following.

 If $f(k)=2$ and $f(k+1)=f(k+2)=1$ we exchange $f$ and $f_2$ on $j',k$ and conclude by Lemma \ref{lem:no02}. If $f(k+2)=2$ and $f(k+1)=f(k)=1$ we exchange $f$ and $f_2$ on $j',k,k+1$. If $f(k)=f(k+2)=0$ and $f(k+1)=1$, we exchange $f$ and $f_2$ on $j',k+2$ and also conclude by Lemma \ref{lem:no02}. If $f(k)=f(k+1)=0$ and $f(k+2)=1$ we can exchange $f$ and $f_3$ on $j',k$ and then $f_2$ and $f_3$ on $k+2,j'$.
 We obtain the following stronger statement.
\wn\label{wn:zg1}
All flows $f\in M_1$ such that $\sum_{i=k}^{k+2}f(i)=1$ agree on $I_1\cup I_3$.$\hfill\square$
\kwn

To finish the proof we have to consider flows $f\in M_1\setminus f_1$ for which $\sum_{i=k}^{k+2}f(i)=0$. By Lemma \ref{lem:no02} and Lemma \ref{lem:more1} we have 
$f(k)=f(k+1)=f(k+2)=1$. 

Now, our strategy is a little different. We do not show directly that each $f$ differs from $h$ on one index. Instead,
by Lemma \ref{lem:more1} we associate to each such $f$ a flow $\tilde f$ that associates only one $1$ to $I_2$. We will show how to modify each pair, so that both flows differ from $h$ on at most one index.
By Lemma \ref{lem:no221} we can further assume that $\tilde f(I_2)=\{0,0,1\}$. By Corollary \ref{wn:zg1} we know that $\tilde f$ agrees with $h$ on $I_1\cup I_3$. 
If $f$ differs from $\tilde f$ on at least two indices from $I_1\cup I_3$, we can exchange $f$ and $\tilde f$ on precisely on index $i\in I_2$ such that $\tilde f(i)=0$ and one or two indices from $I_1\cup I_3$ reducing to previous cases when $\sum_{i=k}^{k+2}f(i)\neq 0$.
This finishes the proof of Proposition \ref{ob:gl} in Case 1 a).

Case 1 b) There is no flow in $M_1$ that associates $0$ to $k+2$ and $1$ to $k$ and $k+1$. In particular, we may assume $f_3(k)=0$ and $f_3(k+1)=1$. 
\lem\label{lem:1agree}
We can assume that all $f\in M_1$ such that $\sum_{i\in I_2}f(i)=1$ agree on $I_1\cup I_3$.
\klem
\dow
It is enough to prove that $f_2$ and $f_3$ agree. Suppose $f_2(i)\neq f_3(i)$. Exchanging $f_2$ and $f_3$ we may assume that $f_2(i)=f_3(i)+1$. Exchanging $f_2$ and $f_3$ on $i, k+2$ we conclude by Lemma \ref{lem:no02}.
\kdow
\lem\label{lem:0agree}
We can assume that all flows $f\in M_1$ such that $\sum_{i\in I_2}f(i)=0$ agree with $h$ on $I_1\cup I_3$, apart from one index.
\klem
\dow
By Lemmas \ref{lem:no02} and \ref{lem:more1} we have $f(I_2)=\{1,1,1\}$. 
Suppose $f(i)\neq f_2(i)$ for $i\in I_1\cup I_3$. If $f(i)=f_2(i)+1$ we can exchange $f$ and $f_2$ on $k,k+1,i$ and conclude that $i$ is the unique index with $f(i)\neq f_2(i)$ by Lemma \ref{lem:1agree}. Suppose $f(i)=f_2(i)+2$. Exchanging $f$ and $f_3$ on $k+2,i$ reduces to Case 1 a).
\kdow
\lem\label{lem:2agree}
We can assume that all flows $f\in M_1$ such that $\sum_{i\in I_2}f(i)=2$ agree with $h$ on $I_1\cup I_3$ apart from one index.
\klem
\dow
By Lemma \ref{lem:no02} the image $f(I_2)$ as a multiset must equal $\{0,1,1\}$ or $\{1,2,2\}$. Suppose it is $\{1,2,2\}$. There is a two element subset $B\subset I_2$ such that $f(B)=\{1,2\}$ and either $f_2(B)=\{0,0\}$ or $f_3(B)=\{0,0\}$. We can exchange $f$ and either $f_2$ or $f_3$ on $B$, obtaining a reduction by Lemma \ref{lem:no02}. Thus we assume $f(I_2)=\{0,1,1\}$. By Case 1 a) we can assume $f(k)=0$ and $f(k+1)=f(k+2)=1$.
Let $i\in I_1\cup I_3$ be such that $f(i)\neq f_2(i)$. If $f(i)=f_2(i)+2$ we can exchange $f$ and $f_2$ on $i, k+1$ obtaining the uniqueness of $i$ by Lemma \ref{lem:1agree} (as after the exchange the flow must exactly agree with $h$ on $I_1\cup I_3$). Thus we assume $f(i)=f_2(i)+1$. If $f_1(i)\neq 0$ we may exchange $f$ and $f_1$ on $i,k$ increasing $k$. Thus we assume $f_1(i)=0$. There must exist one more index $j\in I_1\cup I_3$ such that $f(j)=f_2(j)+1$. By exchanging $f$ and $f_2$ on $k+1,i,j$ by Lemma \ref{lem:1agree} we see that $f$ agrees with $f_2$ on $I_1\cup I_3\setminus\{i,j\}$. As before $f_1(j)=0$. 

By Lemma \ref{lem:f_2zgf_1} we assume that either $f_2(i)=0$ or $f_2(i)=1$. In the latter case $f(i)=2$. We can apply the cubic relation among $f_1,f_2,f$ that on indices $i,k,k+1$ is given by:
$$
\xymatrix{
\txt{{(0,1,1)}\\{(1,0,0)}\\{(2,0,1)}}&=&
\txt{{(2,0,0)}\\{(0,0,1)}\\{(1,1,1)}}
}$$
and all other indices remain unchanged.
This cubic relation increases $k$.

Hence, we may assume $f_2(i)=f_2(j)=0$
. Note that there must exist further indices $i',j'\in I_1$ such that $f_2(i')=f_2(j')=1$. It follows that $f(i')=f(j')=1$. 

Claim: In the situation above, we can assume that each flow in $M_1$ different from $f_1$ attains at least three times $1$ on indices $k+1,k+2,i'j'$.
\dow[Proof of the Claim]
Consider any flow $f'\in M_1$ different from $f_1$.

Suppose $f'(k+1)=f'(k+2)=1$. We have already shown that either $f'$ differs from $f_2$ on $I_1\cup I_3$ on one index, in which case the claim holds, or it disagrees on two indices, none of them equal to $i',j'$, in which case the claim also holds.

Suppose $f'(k+1),f'(k+2)\neq 1$. By Lemma \ref{lem:no02} it either associates twice $0$ or twice $2$. In the latter case we can exchange $f'$ and $f_2$ on $k+1,k+2$ obtaining a contradiction with Lemma \ref{lem:no02}. In the former case, by the same lemma $f'(k)=1$. This gives reduction to Case 1 a).

We are left with the case in which we can assume $f'(k+2)=1\neq f'(k+1)$. We consider two cases.
\begin{enumerate}
\item $f'(k+1)=2$. We have $f'(k)\neq 0$ by Lemma \ref{lem:no02}. If $f'(k)=2$ then we can exchange $f'$ and $f_3$ on $k,k+1$ obtaining contradiction with Lemma \ref{lem:no02}. Hence, $f'(k)=1$. The Claim follows by Lemma \ref{lem:1agree}.
\item $f'(k+1)=0$. If $f'(k)=0$ the Claim follows by Lemma \ref{lem:1agree}. By Lemma \ref{lem:no02} we can assume $f'(k)=1$. This gives a reduction to Case 1 a).
\end{enumerate}
This finishes the proof of the claim.
\kdow
Hence, $f_1$ is the only flow that can attain less than three times $1$ on $k+1,k+2,i',j'$. However, $f$ associates only $1$ to all four indices. Hence, by the Claim the average over all flows in $M_1$ of $1$ on these four indices is at least three. As $g_1\in M_2$ does not have any $1$ on these indices in follows that there must exist a flow $g\in M_2$ that associates $1$ to all four indices. We can exchange $g$ and $g_1$ on $k+1,k+2,i'$ which increases $k$. This finishes the proof.
\kdow
The Lemmas \ref{lem:1agree}, \ref{lem:0agree}, \ref{lem:2agree} finish the proof of Proposition \ref{ob:gl} in Case 1b).

{\bf Case 2)} None of the flows in $M_1$ associates twice $0$ and once $1$ on $I_2$. 
By the case assumption and Lemma \ref{lem:no02} we can assume that:
$$f_2(k)=0,f_2(k+1)=1,f_2(k+2)=1,$$
$$f_3(k)=1,f_3(k+1)=0,f_3(k+2)=1,$$
$$f_4(k)=1,f_4(k+1)=1,f_4(k+2)=0.$$
\lem\label{lem:all}
We can assume that the image of $I_2$  by any flow $f\in M_1$ as a multiset is equal to one of:
$$\{1,1,1\},\{1,1,0\},\{1,1,2\},\{1,2,2\}.$$
\klem
\dow
By Lemma \ref{lem:no02} the only possible multiset that does not contain $1$ is $\{2,2,2\}$, which is excluded by Lemma \ref{lem:more1}. 
If $f(I_2)$ contains exactly one $1$, then by Lemma \ref{lem:no02} it equals either $\{1,2,2\}$ or $\{1,0,0\}$, the latter possibility contradicting Case assumption.
The other possibilities appear in the statement.
\kdow
We will now obtain a contradiction (which proves that we can always reduce to one of the previous cases or increase $k$). In analogy to Lemma \ref{lem:no02}, all flows $g\in M_2$ that do not contain $0$ in the image $g(I_2)$ must attain three times $2$ on $I_2$. To each such flow we can associate a flow $f\in M_1$ such that $f(k+2)=2$. By Lemma \ref{lem:all} all such flows also do not have $0$ in the image of $I_2$. Then we can pair arbitrarily other flows $g\in M_2$, $g\neq g_1$ with flows $f\in M_1$, $f\neq f_1$. By Lemma \ref{lem:all} in each pair $(f,g)$ the flow $g$ attains  $0$ on $I_2$ at least as many times as $f$. Pairing $g_1$ with $f_1$, we obtain a contradiction with compatibility of $M_1$ and $M_2$ counting the number of $0$ on $I_2$.

This finishes the proof of Proposition \ref{ob:gl}.
The proposition assures the existence of a function $h$ for $M_1$ and similarly $h'$ for $M_2$. Notice that on $I_1\cup I_3$ the functions $h$ and $h'$ can disagree only on at most 2 indices, as otherwise, just by comparing the indices on which the functions differ, we would have a contradiction with the fact that $M_1$ and $M_2$ are in a relation. Hence, for $n\geq 7$ we can find two indices $i,j\in I_1\cup I_3$ on which $h$ and $h'$ agree.
\lem\label{lem:two}
For $i,j$ as above the two multisets of pairs:
$$\{(f(i),f(j)):f\in M_1\},\{(g(i),g(j)):g\in M_2\}$$
are equal.
\klem
\dow
By definition the following multisets are equal:
$$\{f(i):f\in M_1\}=\{g(i):g\in M_2\}.$$
If in any of the two multisets an element $o$ different from $h(i)$ appears it gives a pair $(o,h(j))$ in both of the multisets in the statement of the lemma. In the same way we have:
$$\{f(j):f\in M_1\}=\{g(j):g\in M_2\}$$
and each element $o'\neq h(j)$ gives a pair $(h(i),o')$. All the other pairs equal $(h(i),h(j))$.
\kdow
By Lemma \ref{lem:two} we can sum up entries on indices $i,j$ obtaining two compatible multisets $M_1',M_2'$ of flows of size $n-1$. By induction, such relation can be generated and we can lift each quadric and cubic in the generation process. After this procedure, it is enough to apply quadric relations, exchanging only flows on indices $i,j$ to generate the relation represented by $M_1,M_2$.
This finishes the case when $f_1$ attains three times an element different from $0$.

\subsubsection{Let us now assume $k=n-3$, $f_1(k)=f_1(k+1)=1$, $f_1(k+2)=f_1(k+3)=2$.}
\lem\label{exclusions}
We can assume that there is no flow $f'\in M_1$ or $g'\in M_2$ such that:
\begin{enumerate}
\item ($f'(k)=2$ and $f'(k+1)=0$) or ($f'(k)=0$ and $f'(k+1)=2$) or ($f'(k+2)=0$ and $f'(k+3)=1$) or ($f'(k+2)=1$ and $f'(k+3)=0$) or ($f'(i)=0$ and $f'(j)=0$ for $i=k$ or $k+1$ and $j=k+2$ or $k+3$),
\item ($g'(i)=1$ and $g'(j)=2$) where $k\leq i,j$ and either $i=k,k+1$ or $j=k+2,k+3$.
\end{enumerate}
\klem
\dow
As all the statements are similar and easy let us only prove the case $f'(k)=2$ and $f'(k+1)=0$. Then, we can exchange $f'$ and $f_1$ on $k,k+1$.
\kdow
By symmetry we consider two cases.

{\bf Case 1)} There exists a flow $f_2\in M_1$ such that $f_2(k)=f_2(k+1)=0$.
\lem
We may assume $f_2(k+2)=f_2(k+3)=2$.
\klem
\dow
Otherwise we could exchange $f_1$ and $f_2$ and increase $k$.
\kdow
\lem
We can assume that for any flow $f'\in M_1$ such that $f'(k+2)=0$ (resp. $f'(k+3)=0$) we have $f'(k)=f'(k+1)=1$. 
\klem
\dow
By Lemma \ref{exclusions} for $i=k,k+1$ we have $f'(i)=1$ or $f'(i)=2$. If $f'(i)=2$ we can exchange $f'$ and $f_2$ on $i,k+2$ (resp. $k+3$) and conclude by Lemma \ref{exclusions}.
\kdow
Hence, the flows in $M_1$ associate strictly more times $1$ on indices $k,k+1$ than $0$ on indices $k+2, k+3$.
Let $g_2\in M_2$ be a flow that attains strictly more times $1$ on $k,k+1$ than $0$ on $k+2,k+3$, say $g_2(k+1)=1$. By Lemma \ref{exclusions} we have $g_2(k+2),g_2(k+3)\neq 2$. If $g_2(k+2)=g_2(k+3)=1$ we may exchange $g_2$ and $g_1$ on $k+1,k+2,k+3$, increasing $k$. By the choice of $g_2$ we must have $g_2(k)=1$. Notice however that if $g_2(k+2)=1$ (resp. $g_2(k+3)=1$) we can exchange $g_1$ and $g_2$ on $k,k+1,k+2$ (resp. $k+3$), increasing $k$. Thus, $g_2(k+2)=g_2(k+3)=0$, which contradicts the choice of $g_2$.

{\bf Case 2)} There is no flow $f'\in M_1$ or $g'\in M_2$ such that:
\begin{enumerate}
\item $f'(k)=f'(k+1)=0$ or $f'(k+2)=f'(k+3)=0$,
\item $g'(k)=g'(k+1)=1$ or $g'(k+2)=g'(k+3)=2$.
\end{enumerate}
Hence, for any $f'\in M_1$ if $f'(k)=0$ then $f'(k+1)=1$. This contradicts the fact that for any $g'\in M_2$ if $g'(k+1)=1$ then $g'(k)=0$.
\subsection{Suppose $k=n-1$}
As before we let $g_1=\underline 0$. We can also assume $f_1=(0,\dots,0,1,2)$.
The diagram of the proof is as follows.
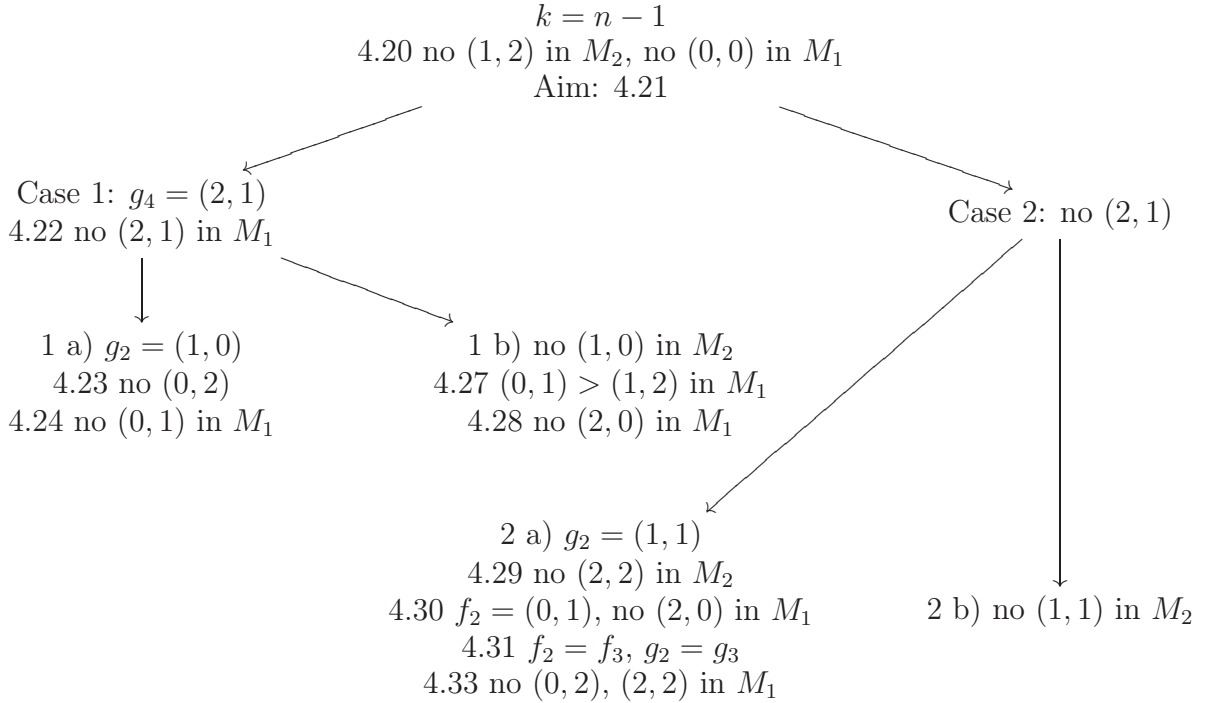
\[
\xymatrix{
&\txt{{$k=n-1$}\\{\ref{lem:ob} no $(1,2)$ in $M_2$, no $(0,0)$ in $M_1$}\\{Aim: \ref{ob:gl2} }}\ar[ld]\ar[rd]\\
\txt{{Case 1: $g_4=(2,1)$}\\{\ref{lem:no21} no $(2,1)$ in $M_1$}}\ar[d]\ar[rd]&&\txt{Case 2: no $(2,1)$}\ar[ldd]\ar[dd]\\
\txt{{1 a) $g_2=(1,0)$}\\{\ref{lem:only22} no $(0,2)$}\\{\ref{lem:2to2} no $(0,1)$ in $M_1$}}&\txt{{1 b) no $(1,0)$ in $M_2$}\\{\ref{lem:more01} $(0,1)>(1,2)$ in $M_1$}\\{\ref{lem:no20M1} no $(2,0)$ in $M_1$}}\\
&\txt{{2 a) $g_2=(1,1)$}\\{\ref{lem:no22} no $(2,2)$ in $M_2$}\\{\ref{lem:no20} $f_2=(0,1)$, no $(2,0)$ in $M_1$}\\{\ref{lem:firstagree} $f_2=f_3$, $g_2=g_3$}\\{\ref{lem:wyklucz} no $(0,2)$, $(2,2)$ in $M_1$}}&\txt{2 b) no $(1,1)$ in $M_2$}\\
}
\]
\dfi[type of a flow]
We say that a flow $f$ is of type $(a,b)$ for $a,b\in \z_3$ if $f(n-1)=a$ and $f(n)=b$.
\kdfi
\lem\label{lem:ob}
If there exists a flow
\begin{enumerate}
\item $g\in M_2$ of type $(1,2)$ or
\item $f\in M_1$ of type $(0,0)$,
\end{enumerate}
then we can increase $k$. We thus assume that such flows do not exist.
\klem
\dow
Follows by obvious quadratic exchange.
\kdow
As $M_1$ and $M_2$ are compatible we can assume $f_2(n-1)=0\neq f_2(n)$, $f_3(n-1)\neq 0=f_3(n)$, $g_2(n-1)=1$, $g_2(n)\neq 2$, $g_3(n)=2$, $g_3(n-1)\neq 1$.
Our aim will be to prove the following proposition.
\ob\label{ob:gl2}
We can assume that 
for each flow $f\in M_1$ different from $f_1$ there exists at most one index $i<n-1$ such that $f(i)\neq f_2(i)$
.
\kob
We consider two cases.

{\bf Case 1)} There exists a flow $g'\in M_1\cup M_2$ of type $(2,1)$. In fact, by the action of $\Gf$ we can assume $g'\in M_2$. We let $g_4=g'$.
\lem\label{lem:no21}
We can assume there is no $f\in M_1$ of type $(2,1)$.
\klem
\dow
In such a case we could exchange the flows so that $f_1(n-1)=g_1(n-1)=2$ and $f_1(n)=g_1(n)=1$ that would make the flows $f_1$ and $g_1$ equal.
\kdow
By Lemma \ref{lem:ob} it is enough to consider two following subcases.

1 a) $g_2(n)=0$
\lem\label{lem:only22}
We can assume that there is no flow $g\in M_2$ of type $(0,2)$. Hence if $g(n)=2$ then $g(n-1)=2$.
\klem
\dow
Assume such $g$ exists. We can make an exhange between $g_4,g,g_2$ that on entries $n-1,n$ is as follows:
$$\xymatrix{\txt{{(2,1)}\\{(0,2)}\\{(1,0)}}&=&\txt{{(1,2)}\\{(2,0)}\\{(0,1)}}}.$$
This, by Lemma \ref{lem:ob} increases $k$. The last statement follows from 
Lemma \ref{lem:ob}.
\kdow
\lem\label{lem:2to2}
We can assume that there is no flow $f\in M_1$ of type $(0,1)$. In particular, $f_2(n)=2$.
\klem
\dow
Suppose such $f$ exists. First notice that we can assume that there is no $\tilde f\in M_1$ such that $\tilde f(n-1)=2$ and $\tilde f(n)=0$. Indeed, in such a case we could make an exchange among $f_1,f,\tilde f$ that on entries $n-1,n$ would be:
$$\xymatrix{\txt{{(1,2)}\\{(0,1)}\\{(2,0)}}&=&\txt{{(2,1)}\\{(1,0)}\\{(0,2)}}},$$
that would allow to increase $k$ by Lemma \ref{lem:no21}.

Hence, also by Lemma \ref{lem:no21} we see that each $f_i\in M_1$ such that $f_i(n-1)=2$ must satisfy $f_i(n)=2$. Taking into account $f_1$ it follows that strictly more flows associate $2$ to $n$ than to $n-1$. This, looking at $M_2$, contradicts Lemma \ref{lem:only22}.
\kdow

\uwa
By Lemma \ref{lem:no21} there must be at least as many flows of type $(2,0)$ as $(1,2)$ in $M_1$. Hence, we may assume that $f_3(n-1)=2$ and $f_2$ and $f_3$ agree on indices $i<n-1$.
\kuwa
We now prove Proposition \ref{ob:gl2} in Case 1 a).
\dow
\emph{Step 1: Flows of type $(1,2)$.}

We will modify the flows of type $(1,2)$ differing from $f_2$ on at least $2$ indices smaller than $n-1$. By Lemma \ref{lem:only22} we know that there are at least as many flows that associate $2$ to $n-1$ as those that associate $2$ to $n$. By Lemma \ref{lem:no21} there must be at least as many flows $\tilde f\in M_1$ such that $\tilde f(n-1)=2$, $\tilde f(n)=0$ as those $f\in M_1$ such that $f(n-1)=1$, $f(n)=2$. Hence, we can pair each $f$ with $\tilde f$. If $f$ differs from $f_2$ (hence also from $f_3$ and $\tilde f$) on some $2$ indices smaller than $n-1$, then by Lemma \ref{lem:zamiana} we can exchange $f$ and $\tilde f$ on $n$, keeping $n-1$.

\emph{Step 2: Flows such that $f(n-1)=0$ or $f(n)=0$.}

Here, the conclusion of the lemma follows from Lemma \ref{lem:zamiana} by exchanging with $f_2$ or $f_3$. By previous lemmas the only cases left are flows of type $(1,1)$ or $(2,2)$. 

\emph{Step 3: Flows of type $(1,1)$.}

We can assume that each flow of type $(1,1)$ agrees with $f_2$ on all indices $i<n-1$. Indeed, otherwise, it it would differ on at least two indices, which we could exchange, obtaining at least two disagreements between $f_2$ and $f_3$.  By Lemma \ref{lem:zamiana} we could then obtain the flow of type $(0,0)$.
Moreover, at least one flow $f_j\in M_1$ of type $(1,1)$ must exist as only such flows satisfy $f_j(n)=1$. 

\emph{Step 4: Flows of type $(2,2)$.}

If a flow of type $(2,2)$ differs from a flow of type $(1,1)$ (which is equivalent to differing from $f_2$) on at least two indices $i<n-1$ we would be able to make an exchange by Lemma \ref{lem:zamiana} obtaining a flow of type $(2,1)$. This finishes the proof by Lemma \ref{lem:ob}.
\kdow

1 b) There is no flow $g\in M_2$ such that $g(n-1)=1$ and $g(n)=0$. In particular, $g_2(n)=1$. Notice that if $g\in M_2$ and $g(n-1)=1$, then $g(n)=1$.
\wn\label{cor1}
Strictly  more flows in $M_1$ associate $1$ to $n$ than to $n-1$. $\hfill\square$
\kwn
\lem\label{lem:more01}
There are strictly more flows $f\in M_1$ of type $(0,1)$, than those of type $(1,2)$. In particular, we can assume $f_2(n)=1$.
\klem
\dow
Follows from Corollary \ref{cor1} by Lemma \ref{lem:no21}.
\kdow

\lem\label{lem:no20M1}
We can assume that there are no flows $f\in M_1$ of type $(2,0)$. In particular, $f_3(n-1)=1$.
\klem
\dow
Indeed, in such a case we can exchange $f, f_2,f_1$ on indices $n-1,n$ as follows:
$$
\xymatrix{
\txt{{(2,0)}\\{(0,1)}\\{(1,2)}}&=&\txt{{(0,2)}\\{(1,0)}\\{(2,1)}}}$$
and conclude by Lemma \ref{lem:no21}.
\kdow
We now prove Proposition \ref{ob:gl2} in Case 1 b).
\dow
First notice that flows of type either $(0,1)$ or $(1,0)$ or $(2,2)$ must agree on indices smaller than $n-1$. Indeed, if such flows differ, they have to differ on at least two indices, that would allow to make an arbitrary exchange on indices $n-1,n$ by Lemma \ref{lem:zamiana}. There is always an exchange that allows to increase $k$.

We now consider flows $f\in M_1$, different from $f_1$, of type $(1,2)$. By Lemma \ref{lem:more01} we can pair each such $f$ with $\tilde f$, of type $(0,1)$. If $\tilde f$ and $f$ differ on at least $2$ indices smaller than $n-1$, by Lemma \ref{lem:zamiana} we can make an exchange that on indices $n-1,n$ is as follows:
$$\xymatrix{\txt{{(0,1)}\\{(1,2)}}&=&\txt{{(1,1)}\\{(0,2)}}}.$$

The only flows left to consider are of type $(1,1)$ or $(0,2)$. Notice that there must exist a flow $f_j\in M_1$ of type $(2,2)$, as only such flow satisfies $f_j(n-1)=2$. If a flow of type $(1,1)$ exists and differs on at least two indices smaller than $n-1$ we can make an exchange by Lemma \ref{lem:zamiana} with $f_j$:
$$\xymatrix{\txt{{(1,1)}\\{(2,2)}}&=&\txt{{(1,2)}\\{(2,1)}}},$$
that allows to increase $k$. If a flow $f'\in M_1$ of type $(0,2)$ exists and differs on at least two indices smaller than $n-1$ then we can find a subset $I\subset [0,n-2]$ such that $\sum_{i\in I}f'(i)=\sum_{i\in I}f_2(i)+1$. We can then make an exchange among $f',f_3,f_1$ that on indices $I,n-1,n$ is as follows:
$$\xymatrix{\txt{{(a+1,0,2)}\\{(a,1,0)}\\{(\underline 0,1,2)}}&=&\txt{{(a,1,2)}\\{(a+1,1,2)}\\{(\underline 0,0,0)}}},$$
where we represent all values of $f_3$ on $I$ by $a$.
\kdow
{\bf Case 2)} There is no flow $f$ among $M_1\cup M_2$ such that $f(n-1)=2$ and $f(n)=1$.
We consider two subcases.

Case 2 a) There exists a flow $g\in M_2$ such that $g(n-1)=g(n)=1$. In this case we let $g_2=g$.
\lem\label{lem:no22}
If there is a flow $g\in M_2$ of type $(2,2)$ then we can increase $k$. In particular we may assume $g_3(n-1)=0$.
\klem
\dow
In such a case we can apply the relation among $g_1,g_2,g$ that on indices $n-1,n$ is of the following form:
$$\xymatrix{\txt{{(0,0)}\\{(1,1)}\\{(2,2)}}&=&\txt{{(1,2)}\\{(2,0)}\\{(0,1)}}}.$$
\kdow
\lem\label{lem:no20}
We may assume that $f_2(n)=1$.
Moreover, if there is a flow $f\in M_1$ of type $(2,0)$ then we can make a reduction. In particular, we may assume $f_3(n-1)=1$.
\klem
\dow
First let us note that the first sentence implies the second. Indeed, we can make an exchange among $f_1,f,f_2$ that on indices $n-1,n$ would be:
$$\xymatrix{\txt{{(1,2)}\\{(2,0)}\\{(0,1)}}&=&\txt{{(0,0)}\\{(1,1)}\\{(2,2)}}}.$$

To prove the first sentence assume the contrary that for any $f'\in M_1$ if $f'(n-1)=0$ then $f'(n)=2$. It follows that strictly more flows associate $2$ to $n$ than $0$ to $n-1$. This contradicts Lemma \ref{lem:no22}.
\kdow
\lem\label{lem:firstagree}
We can assume that:
\begin{enumerate}
\item the flows $f_2$ and $f_3$ agree on all indices smaller than $n-1$ and
\item the flows $g_2$ and $g_3$ agree on all indices smaller than $n-1$.
\end{enumerate}
\klem
\dow
Each pair, if it differs, it must differ on at least $2$ indices and the conclusion follows from Lemma \ref{lem:zamiana}.
\kdow
\lem\label{lem:11problem}
We can assume that each flow $f\in M_1$ different from $f_1$ either:
\begin{enumerate}
\item differs from $f_2$ only on one index smaller than $n-1$, or
\item $f(n-1)=f(n)=1$ and $f$ differs from $f_2$ on exactly two indices $i,j<n-1$ for which $f(i)=f_2(i)+1$, $f(j)=f_2(j)+1$
\end{enumerate}
\klem
\dow
We only have to consider flows of types:
$(0,1)$, $(0,2)$, $(1,0)$, $(1,1)$, $(1,2)$, $(2,2)$. From Lemma \ref{lem:firstagree} we obtain that $(0,1)$, $(1,0)$ and $(2,2)$ agree on all entries smaller than $n-1$. If a $(0,2)$ flow differs on two indices, we can increase $k$ by Lemma \ref{lem:zamiana}. If any flow $f$ differs on at least $3$ indices we can exchange $f,f_2,f_3$ obtaining a flow that assigns $0$ both to $n-1$ and $n$. Hence, to prove the Lemma, we only need to consider flows $f$ of type $(1,2)$ for which there are two indices $i',j'<n-1$ such that $f(i')=f_2(i')+2$, $f(j')=f_2(j')+2$.

Flows of type $(1,2)$ can be paired with $(0,1)$ flows (if there are not enough $(0,1)$ flows we can apply Lemma \ref{lem:no20}) and in each pair we can apply the relation:
$$\xymatrix{\txt{{(a,0,1)}\\{(a+2,1,2)}}&=&\txt{{(a+2,1,1)}\\{(a,0,2)}}},$$
where $a=f_2(i')$.
\kdow
\lem\label{lem:wyklucz}
Either all flows in $M_1$ apart from $f_1$ differ from $f_2$ on at most one index smaller than $n-1$ or we can assume there are no flows $f\in M_1$ of type $(0,2)$ or $(2,2)$.
\klem
\dow
Suppose there exists a flow $f$ described in point 2) of Lemma \ref{lem:11problem}.
If a flow of type $(0,2)$ exists we can exchange it with $f$ on last two entries obtaining a contradiction with Lemma \ref{lem:11problem}.

If a flow of type $(2,2)$ exists we can make an exchange obtaining a flow of type $(2,1)$, hence reduce to Case 1).
\kdow
To finish the proof of Proposition \ref{ob:gl2} 
we have to prove that there exists a flow in $M_1$ of type $(0,2)$ or $(2,2)$. Suppose this is not the case - we will reach a contradiction.
All flows in $M_1$ are of one of the types $(1,2)$, $(0,1)$, $(1,0)$, $(1,1)$. It follows that the flows in $M_2$ can be only of types $(0,0)$, $(1,1)$, $(0,2)$, $(1,0)$, $(0,1)$.

Choose an index $i<n-1$ such that $f_2(i)\neq g_2(i)$. Suppose there are $x$ flows of type $(1,2)$ in $M_1$. Then there are also $x$ flows of type $(0,2)$ among $M_2$. Let there be $y$ flows of type $(0,1)$ and $z$ of type $(0,0)$ in $M_2$. Then there are $x+y+z$ flows of type $(0,1)$ in $M_1$. Let $w$ be the number of flows of type $(1,1)$ in $M_1$. Then there are $x+z+w$ flows of type $(1,1)$ in $M_2$. Suppose there are $q$ flows of type $(1,0)$ in $M_2$. Then there are $z+q$ flows of type $(1,0)$ in $M_1$. Hence, on index $i$, $f_2(i)$ appears at least $x+z+y+z+q$ times and $g_2(i)$ appears at least $x+z+w+x$ times. However, the sum of these two numbers $3x+3z+y+q+w$ is larger than the cardinality of the multisets $2z+2x+w+q+y$, which gives the contradiction.

Case 2 b) There is no flow $g\in M_2$ such that $g(n-1)=g(n)=1$. By the action of $\Gf$ we can reduce to Case 2 a) also when:
\begin{enumerate}
\item there is a flow of type $(0,1)$ in $M_1$ (by the nontrivial automorphism of $\z_3$),
\item there is a flow of type $(2,0)$ in $M_1$ (by the transposition of $n-1$ and $n$).
\end{enumerate}
Recall that we also know that there are no flows of type $(0,0)$ or $(2,1)$ in $M_1$. 
By the case assumption we know that in $M_2$ there are no flows of type $(1,1)$, $(1,2)$. Hence, strictly more flows in $M_2$ associate $0$ to $n$ than $1$ to $n-1$. This contradicts the fact that there are no flows of type $(2,0)$ and $(0,0)$ in $M_1$.

This finishes the proof of Proposition \ref{ob:gl2}. 

Among all the multisets that can be reached from $M_2\setminus g_1$ by exchanging at most three flows, let $\tilde M_2$ be one (possibly out of many) that minimizes the number $s$ of indices smaller than $n-1$ on which a flow from $\tilde M_2$ can differ from $f_2$. Further, it also minimizes the number of flows differing on exactly $s$ indices, and further on $s-1$ indices.
First note that $s\leq 3$. Indeed, if there is a flow $g\in \tilde M_2$ differing on $4$ indices, by compatibility we can find a flow $\tilde g\in \tilde M_2$ that equals $f_2$ on $[n-2]$ and we can make an exchange on $2$ or $3$ indices. 

Suppose $s=3$. By compatibility we can find a flow $\tilde g\in \tilde M_2$ that equals $f_2$ on $[n-2]$ (there exist at least three such flows). Let $\tilde g$ be of type $(a,b)$.
Without loss of generality we can assume that $g\in \tilde M_2$ and $g(l)=f_2(l)+1$ for $l=1,2,3$. There can be no flows in $\tilde M_2$ that attain $f_2(i)+2$ at index $i$ for any $i<n-1$, as we would be able to make an exchange with $g$ and $\tilde g$.
By the definition of $\tilde M_2$ all flows in that multiset that differ from $f_2$ on $3$ indices in $[n-2]$ must also be of type $(a,b)$. So must all the flows that agree with $f_2$. Further, flows that differ on $1$ or $2$ indices must either have $a$ on $n-1$ or $b$ on $n$. By compatibility there must exist a flow in $M_1$ of type $(a,b)$. We may assume that $f_2$ is of type $(a,b)$ and equals $\tilde g$.

Suppose $s=2$. Without loss of generality there are at least as many flows in $\tilde M_2$ that differ from $f_2$ on indices $1,2$ as on indices $4,5$. We pair each flow that differs on $4,5$ with a flow that differs on $1,2$. By Lemma \ref{lem:zamiana} we can make an exchange in each pair, in such a way that none of the flows differs both on $4$ and $5$. Then for $i,j=4,5$ we proceed exactly as in Lemma \ref{lem:two} and the paragraph after it. 

If $s=1$ or $s=0$ we may proceed as above, choosing any $i,j\leq n-2$. This finishes the proof of the main theorem.

\bibliographystyle{amsalpha}
\bibliography{Xbib}
\noindent Mateusz Micha{\l}ek,
Polish Academy of Sciences, Warsaw, Poland,
{\tt wajcha2@poczta.onet.pl}
\end{document}